\documentclass[12pt,twoside]{amsart}
\usepackage{amssymb}

\newtheorem{Thm}{Theorem}
\newtheorem{Cor}{Corollary}
\newtheorem{Lemma}{Lemma}
\newtheorem{Prop}{Proposition}

\theoremstyle{definition}
\newtheorem{Defn}{Definition}
\newtheorem{Ex}{Example}

\theoremstyle{remark}
\newtheorem{Remark}{Remark}

\newcommand{\norm}[1]{\left\Vert#1\right\Vert}
\newcommand{\abs}[1]{\left\vert#1\right\vert}
\newcommand{\chf}{\ensuremath{\mathbf{1}}}
\DeclareMathOperator{\Var}{Var}
\newcommand{\set}[1]{\left\{#1\right\}}
\newcommand{\eps}{\varepsilon}

\renewcommand{\phi}{\varphi}

\newcommand{\NC}{\mathit{NC}}
\newcommand{\Res}{\mathcal{R}}
\newcommand{\m}{{\mathfrak{m}}}

\DeclareMathOperator{\Rtr}{\mathrm{R}}

\newcommand{\mf}[1]{\mathbb{#1}}
\newcommand{\mc}[1]{\mathcal{#1}}

\title{It\^o formula for free stochastic integrals}
\author[M. Anshelevich]{Michael Anshelevich}
\thanks{This work was supported in part by an NSF postdoctoral fellowship}
\address{Department of Mathematics, University at California, Berkeley, CA 94720, USA}
\email{manshel@math.berkeley.edu}
\subjclass{Primary 46L54; Secondary 60G10, 81S25}

\date{February 7, 2001}

\begin{document}

\begin{abstract}
The objects under investigation are the stochastic integrals with respect to free Levy processes. We define such integrals for square-integrable integrands, as well as for a certain general class of bounded integrands.  Using the product form of the It\^o formula, we prove the full functional It\^o formula in this context. 
\end{abstract}

\maketitle

\section{Introduction}

The subject of this paper is stochastic integration in the context of free probability. Noncommutative stochastic processes with freely independent increments, especially the free Brownian motion, have been investigated in a number of sources, see \cite{BS98}, \cite{Ans00} and their references. In the latter paper, we started the analysis of such processes, which we also call free stochastic measures, using the combinatorial machinery inspired by the work of Rota and Wallstrom \cite{RW97}. Starting with a free stochastic measure $X(t)$, in that paper we defined a family of multi-dimensional free stochastic measures derived from it, indexed by set partitions. In particular, we defined the family of higher diagonal stochastic measures $\Delta_k(t)$, which give a precise meaning to the heuristic expression $d \Delta_k = (dX)^k$.

In this paper we define integrals with respect to free stochastic measures, and investigate their properties. We restrict the analysis to free stochastic measures consisting of bounded operators. Note that this has no analog in the classical stochastic integration theory: there are no (non-trivial) compactly supported infinitely divisible distributions, while the class of compactly supported freely infinitely divisible distributions is dense, and includes the free analogs of the normal and the Poisson distribution. After completing this paper, we learned from Steen Thornbj{\o}rnsen about a recent preprint \cite{BNT00}. In it, using a remarkable bijection between the free and classical infinitely divisible distributions, the authors define stochastic integrals with respect to any stationary stochastic process with free increments, with the Riemann sums defining the integrals converging in probability. On the other hand, since we are dealing with bounded operators, we can achieve convergence in the operator norm, both for the integrals and the limits defining the higher diagonal measures. This requires a definition of a family of ``mixed-$p$'' norms on the integrands, with the integrands which are bounded in the $\infty$-norm giving stochastic integrals which are bounded in the operator norm. Boundedness of the integrators also allows us to define the integrals for a significantly wider (not necessarily scalar, not necessarily continuous) class of integrands.

Most importantly, we prove the functional It\^{o} formula for such integrals, which involves the integration with respect to the original free stochastic measure as well as its higher diagonal measures. The importance of free It\^{o} formulas is indicated by recent applications of the It\^o formula for the free Brownian motion to fine properties of random Gaussian matrices in \cite{CD99, CDG00}. The paper is organized as follows. After some preliminaries in Section \ref{sec:prel}, we define integration with respect to free stochastic measures in Section \ref{sec:int}. Using the It\^{o} isometry, the integration is defined for integrands which are square-integrable in the appropriate sense. Section \ref{sec:ito} contains the various It\^o formulas, which are the main results of the paper.  In Section \ref{sec:linf} we define the ``mixed-$p$''-norms on the integrands, and extend the definition of the integral and the It\^{o} formulas to this context.

\section{Preliminaries}
\label{sec:prel}
In this section we collect the definitions and results from \cite{Ans00} and \cite{BS98}. For consistency, whenever possible we follow their notation.

\subsection{Probability space} 
All the work will proceed in a $W^\ast$-noncommutative probability space $(\mc{A}, \phi)$. Here, $\mc{A}$ is a finite von Neumann algebra, and $\phi$ is a faithful normal trace state on it. As usual, we define the $p$-norms on $\mc{A}$ by $\norm{S}_p = \phi [\abs{S}^p]^{1/p}$, $\norm{S}_\infty = \norm{S}$, and $L^p(\mc{A}, \phi)$ are the completions of $\mc{A}$ with respect to these norms. We say that a self-adjoint element $Y \in \mc{A}$ has distribution $\nu$ (with respect to the trace $\phi$) if $\phi[Y^n] = \int s^n d\nu(s)$ for $n \in \mf{N}$. Note that for bounded $Y$ its distribution has compact support.

Let $\mc{A}$ be filtered by an increasing family of von Neumann subalgebras $\set{\mc{A}_t}_{t \in [0, \infty)}$.

\subsection{$R$-transforms}
Our main references for the background in free probability are \cite{VDN92} and \cite{Voi00b}; see also the references in \cite{Ans00}. Let $\nu$ be a compactly supported probability distribution on $\mf{R}$. Let $G_{\nu}(z) = \int_{\mf{R}} \frac{d \nu(s)}{z-s}$ be its Cauchy transform. It has a power series expansion around infinity 
\begin{equation}
\label{mom}
G_{\nu}(z) = \sum_{k=0}^\infty m_k(\nu) z^{-(k+1)},
\end{equation}
where $m_k(\nu)$ will always denote the $k$-th moment of $\nu$. Define the formal power series $\Rtr_{\nu}(z)$, the \emph{$R$-transform} of $\nu$, by the equation 
\begin{equation}
\label{trans}
G_\nu \left( \frac{1}{z} + \Rtr_{\nu}(z) \right) = z.
\end{equation}
Then $\Rtr_\nu$ is can be extended to a function which is analytic in a Stolz angle $\Gamma$ at $0$ in $\mf{C}^+$ and maps $\Gamma$ into $\mf{C}^+$, and also satisfies $\Rtr_\nu(\bar{z}) = \overline{\Rtr_\nu(z)}$.

Call the coefficients $\set{r_i(\nu)}_{i=1}^\infty$ in the expansion 
\begin{equation}
\label{cum}
\Rtr_\nu(z) = \sum_{k=1}^\infty r_k(\nu) z^{k-1}
\end{equation}
the \emph{free cumulants} of $\nu$. Substituting the equations \eqref{mom} and \eqref{cum} into the functional relation \eqref{trans} and equating the coefficients, we obtain the following relation between the moments and the free cumulants. Let $\NC(n)$ be the lattice of noncrossing partitions of the set of $n$ elements. Let $\pi \in \NC(n)$, and denote $R_\pi(\nu) = \prod_{B \in \pi} r_{\abs{B}}(\nu)$. Then $m_n(\nu) = \sum_{\pi \in \NC(n)} R_{\pi}(\nu)$. Finally, denote $\norm{\nu} := \sup \set{\abs{x}: x \in supp(\nu)}$. Note that $\norm{\nu} = \norm{Y}$ if $Y$ has distribution $\nu$.

Let $\mu$ be a freely infinitely divisible
distribution with compact support. It has the property that $\Rtr_\mu$ can be extended, with the above properties, to the whole $\mf{C}^+$. Normalize $\mu$ so that $\Var (\mu) = 1$. Let $\set{\mu_t}_{t \in [0, \infty)}$ be the corresponding additive free convolution semigroup. It is  determined by the requirement that $\Rtr_{\mu_t}(z) = t \Rtr_{\mu}(z)$. In particular, the $i$'th free cumulant of $\mu_t$ is $t r_i(\mu)$. Note that by assumption $r_2 = 1$. In the sequel, we fix $\mu$ and $\set{\mu_t}_{t \in [0, \infty)}$, and denote the moments of $\mu_t$ by $m_n(t)$, the free cumulants of $\mu$ by $r_n$, and the $R$-transform of $\mu$ by $\Rtr$.

\begin{Defn}
\label{defn:fip}
A \emph{noncommutative stochastic measure} is a map from the
set of finite half-open intervals $I = [a, b) \subset [0, \infty)$ to the self-adjoint part of
$(\mc{A}, \phi)$,
(which can be extended to the map on all Borel subsets)
$I \mapsto X(I)$.
\end{Defn}
The terms ``stochastic measure'' and ``stochastic process'' will be used interchangeably.
We put the following 3 requirements on the stochastic measure:
\begin{enumerate}
\item Additivity: $I_1 \cap I_2 = \emptyset, I_1 \cup I_2 = J \Rightarrow X(I_1) + X(I_2) = X(J)$,
\item Stationarity: The distribution of $X(I)$ depends only on $\abs{I}$, and equals $\mu_{\abs{I}}$,
\item Free increments: for a family $\set{I_1, I_2, \ldots, I_n}$ of disjoint intervals, the corresponding family of operators $\set{X(I_1), X(I_2), \ldots, X(I_n)}$ are freely independent.
\end{enumerate}
We denote $X(t) = X([0, t))$, and say that $X$ is a \emph{free stochastic measure with distribution $\mu$}. Note that $X(0) = 0$.

\begin{Lemma}
\label{lem:norm}
Let $X$ be a free stochastic measure with distribution $\mu$.

\begin{enumerate}
\item The function $t \mapsto X(t)$ is continuous in $L^p(\mc{A}, \phi)$ for all $p < \infty$. 
\item Suppose all the free cumulants of $\mu$ are non-negative (see Section \ref{sec:linf} for further use of this assumption). Under this condition, the function $t \mapsto X(t)$ is not continuous in $L^\infty(\mc{A}, \phi)$ unless $X$ is the free Brownian motion.
\end{enumerate}
\end{Lemma}
\begin{proof}
For part (1), by stationarity, it suffices to prove that $\lim_{t \rightarrow 0} \norm{X_t}_p = 0$. Since the result is to be proven for all $p$, it suffices to do so for all even integer $p$. Moreover, the statement is not about the stochastic measure but about its distributions $\mu_t$. Namely, what needs to be shown is that $\lim_{t \rightarrow 0} \sqrt[p]{m_p(t)} = 0$. Indeed, 
\begin{equation*}
m_p(t) = \sum_{\pi \in \NC(p)} R_\pi t^{\abs{\pi}} = t r_p + o(t).
\end{equation*}
Thus $\sqrt[p]{m_p(t)} = \sqrt[p]{r_p} \sqrt[p]{t} + o(\sqrt[p]{t}) \rightarrow 0$ as $t \rightarrow 0$. 

Under the positivity assumption on the free cumulants, 
\begin{equation*}
\limsup_{t \rightarrow 0^+} \norm{X_t} = \limsup_{t \rightarrow 0^+} \lim_{p \rightarrow \infty} \norm{X_t}_p \geq \limsup_{p \rightarrow \infty} \sqrt[p]{\abs{r_p}}.
\end{equation*}
This is $0$ iff the radius of convergence of the power series defining $\Rtr$ is infinite, in other words if $\Rtr$ is analytic in the complex plane. But the only functions that are entire, map $\mf{C}^+$ into $\mf{C}^+$ and satisfy $F(\bar{z}) = \overline{F(z)}$ are of the form $F(z) = a z + b$, and these are the $R$-transforms of non-scaled, non-centered semicircular distributions.
\end{proof}

\begin{Defn}
A stochastic measure  is adapted to the filtration $\set{\mc{A}_t}_{t \in [0, \infty)}$ of $\mc{A}$ if
\begin{enumerate}
\item $I \subset [0,t) \Rightarrow X(I) \in \mc{A}_t$,
\item $I \cap [0, t) = \emptyset \Rightarrow \mc{A}_t$ and $X(I)$ are freely independent.
\end{enumerate}
\end{Defn}
From now on, we fix an adapted free stochastic measure $X$ with distribution $\mu$.

\begin{Defn}
\label{defn:pr}
Let $I = [a,b) \subset [0, \infty)$. Let $X_i^{(N)} = X([a + \frac{i-1}{N}(b-a), a + \frac{i}{N}(b-a)))$, for $1 \leq N$, $1 \leq i \leq N$. Then we define the $k$-th diagonal measure of $X$ by 
\begin{equation*}
\Delta_k (I) = \lim_{N \rightarrow \infty}
\sum_{i=1}^N (X_i^{(N)})^k,
\end{equation*}
where the limit is in the operator norm. The limit exists by the results in Section 6 of \cite{Ans00}.
\end{Defn}
As before, denote $\Delta_k(t) = \Delta_k([0, t))$.

\begin{Lemma}
If $X$ is an adapted free stochastic measure then so is $\Delta_k$. 
\end{Lemma}
\begin{proof}
Corollary 3 of \cite{Ans01b}.
\end{proof}
The distribution of $\Delta_k$ is determined by the following result:

\begin{Lemma}
\label{lem:r_n}
The free cumulants of the $k$-th diagonal measure of the
process are given by
\begin{equation*}
r_n(\Delta_k(t)) = t r_{nk}.
\end{equation*}
In particular, $\phi[\Delta_k(t)] = t r_k$. 
\end{Lemma}
\begin{proof}
Theorem 2 of \cite{Ans00}.
\end{proof}

\section{Stochastic integrals}
\label{sec:int}
The definition of stochastic integrals proceeds in the usual way: first we define it for simple processes, and then extend it to $L^2$ processes using the It\^{o} isometry. Actually, the objects we will be integrating are biprocesses.

\subsection{Biprocesses}
A biprocess is a map $[0, \infty) \rightarrow \mc{A} \otimes \mc{A}^{op}$. Here $\mc{A}^{op}$ is the opposite algebra of $\mc{A}$. A simple biprocess is piecewise constant and is $0$ outside of a finite interval. Equivalently, simple biprocesses are spanned by all $A \otimes B \chf_{[a, b)}$ for $A, B \in \mc{A}$ and $\chf_{[a, b)}$ the characteristic function of the interval. A biprocess $U$ is adapted if $U(t) \in \mc{A}_t \otimes \mc{A}_t^{op}$ for all $t$. Denote the space of all simple adapted biprocesses by $\mc{B}$; it is spanned by all $A \otimes B \chf_{[a, b)}$ as above with $A, B \in \mc{A}_a$. Throughout, we will only consider the integration of adapted biprocesses. Note that the property of being adapted has a stronger meaning for a stochastic measure than for a biprocess.

Let $U \in \mc{A} \otimes \mc{A}^{op}$. Define the involution on $\mc{A} \otimes \mc{A}^{op}$ by its action on elementary tensors: $(A \otimes B)^\ast = B^\ast \otimes A^\ast$, and extend by linearity. Denote by $\norm{U}_\infty = \norm{U}$ the operator norm of $U$. Note that this is not a $C^\ast$-norm for the above involution, but $\mc{A} \otimes \mc{A}^{op}$ is still a Banach $^\ast$-algebra. Denote by $\langle \cdot, \cdot \rangle$ the inner product in $L^2(\mc{A}, \phi) \otimes L^2(\mc{A}^{op}, \phi)$, given on elementary tensors by $\langle A \otimes B, C \otimes D \rangle = \phi[C^\ast A] \phi[B D^\ast]$, or using the notation below $\langle V, U \rangle = (\phi \otimes \phi) [U^\ast \sharp V]$. Denote by $\norm{\cdot}_2$ the corresponding norm on $\mc{A} \otimes \mc{A}^{op}$.

For a biprocess, we will use the notation $\norm{\norm{U}_p}_q$ to denote $\sqrt[q]{\int_0^\infty \norm{U(s)}_p^q ds}$, including $\norm{\norm{U}_p}_\infty = \sup_{s \in [0, \infty)} \norm{U(s)}_p$. Denote by $\mc{B}^p$ the completion of $\mc{B}$ with respect to the norm $\norm{\norm{\cdot}_\infty}_p$ and by ${\mc{B}^2}'$ the completion of $\mc{B}$ with respect to the norm $\norm{\norm{\cdot}_2}_2$. Note that $\mc{B}^\infty$ does not consist of all bounded adapted biprocesses, but only of those that converge to $0$ (in norm) as time goes to infinity.

\subsection{Notation}
\begin{enumerate}

\item Let $\sharp$ denote the left bimodule action $\sharp: (\mc{A} \otimes \mc{A}^{op}) \times \mc{A} \rightarrow \mc{A}$. Also, for $k \geq 1, 1 \leq i \leq k$ define the right bimodule action $\sharp_{k, i}: \mc{A}^{k+1} \times (\mc{A} \otimes \mc{A}^{op}) \rightarrow \mc{A}^{k+1}$ as follows. For $U \in \mc{A} \otimes \mc{A}^{op}$, $\cdot \sharp_{k, i} U$ is the right multiplication by $1^{\otimes(i-1)} \otimes U \otimes 1^{\otimes(k + 1 - i)}$ in $\mc{A}^{\otimes i} \otimes \mc{A}^{op} \otimes \mc{A}^{\otimes(k-i)}$. Note that for a fixed $k$ and different $i$, the actions $\sharp_{k,i}$ commute. 
Denote $\sharp \m_k: \mc{A}^{\otimes (k+1)} \times \left(\mc{A} \otimes \mc{A} \right)^{\otimes k} \rightarrow \mc{A}^{\otimes (k+1)}$ the right action
\begin{equation*}
\sharp \m_k(U_1, U_2, \ldots, U_k) = \sharp_{k,1} U_1 \sharp_{k,2} U_2 \cdots \sharp_{k,k} U_k.
\end{equation*}
For convenience, write $\m_1(U)$ simply as $U$. Note that the definitions of $U \sharp V$ are consistent.

\item For $k \geq 2$, define the contractions
$\phi_k: \mc{A}^{\otimes k} \rightarrow \mc{A} \otimes \mc{A}$ by $\phi_k = I \otimes \phi^{\otimes (k-2)} \otimes I$. 

\item Finally, denote
$U \otimes_2 V = \phi_3 [1^{\otimes3} \sharp \m_2(U,V)] = \phi_3[(U \otimes 1)(1 \otimes V)]$.

\end{enumerate}

\subsection{Free difference quotient}
Let $C[x]$ be the algebra of complex polynomials in an indeterminate $x$. Define the canonical derivation $\partial: C[x]  \rightarrow C[x] \otimes C[x]$ by the requirement that $\partial(1) = 0$, $\partial(x) = 1 \otimes 1$. Moreover, define the maps $\partial^k: C[x] \rightarrow C[x]^{\otimes(k+1)}$ by $\partial^k = k (1 \otimes \ldots \otimes 1 \otimes \partial) \partial^{k-1}$. More explicitly, on monomials their action is 
\begin{equation*}
\partial^k x^n = k! \sum_{i(0), i(1), \ldots, i(k) = 0}^{i(0) + i(1) + \ldots + i(k) = n-k} x^{i(0)} \otimes x^{i(1)} \otimes \cdots \otimes x^{i(k)}.
\end{equation*}

Derivation property of $\partial$ implies the following property of $\partial^k$: for any polynomial $p$, 
\begin{align*}
\partial^k (p(x) x) &= k (1 \otimes \ldots \otimes 1 \otimes \partial) (\partial^{k-1}(p(x) x)) \\
&= k (1 \otimes \ldots \otimes 1 \otimes \partial) (\partial^{k-1}(p(x)) (1 \otimes \ldots \otimes 1 \otimes x) + (k-1) \partial^{k-2}(p(x)) \otimes 1) \\
&= (1 \otimes \ldots \otimes 1 \otimes x) \partial^k(p(x)) + k \partial^{k-1}(p(x)) \otimes 1,
\end{align*}
where we have used induction in the intermediate step.
In what follows, whenever no confusion arises we will frequently abuse notation by writing $\partial^k(p)(M)$ or $\partial^k (p(M))$ in place of $\partial^k(p)(M, M, \ldots, M)$.

\begin{Defn}
For an adapted free stochastic measure $X$ and an adapted biprocess $U$, we define the stochastic integral $\int U(s) \sharp dX(s)$ by the rule that if $U = A \otimes B \chf_{[a, b)}$ then  
\begin{equation*}
\int U(s) \sharp dX(s) = A (X(b) - X(a))B,
\end{equation*}
and extend by linearity. Denote $\int_{t_1}^{t_2} U(s) \sharp dX(s) = \int (U(s) \chf_{[t_1, t_2)}(s)) \sharp dX(s)$. 
\end{Defn}
Note that $(\int U(s) \sharp dX(s))^\ast = \int U^\ast(s) \sharp dX(s)$. Also, we will often omit the variable of integration and write simply $\int U \sharp dX$ for $\int U(s) \sharp dX(s)$. Whenever we denote by $M$ the integral $\int U(s) \sharp dX(s)$ we will denote by $M(t)$ the integral $\int_0^t U(s) \sharp dX(s)$. Clearly the process $M(t)$ is adapted.

\subsection{Extension to $L^2$}
\label{sec:l2}
In this section we extend the stochastic integral map to adapted biprocesses which are square-integrable in a particular sense.

\begin{Defn}
Let $V, U$ be adapted biprocesses. For $a \in \mf{R}$, define an inner product
\begin{equation*}
\langle V, U \rangle_a' = \int \langle V(s), U(s) \rangle ds + a^2 \langle \int V(s) ds, \int U(s) ds \rangle.
\end{equation*}
Denote by $\norm{U}_{2,a}' = \sqrt{\norm{\norm{U}_2}_2^2 + a^2 \norm{\int U(s) ds}_2^2}$ the corresponding norm, and by ${\mc{B}^{2, a}}'$ the completion of $\mc{B}$ with respect to this norm. 
\end{Defn}

\begin{Prop}
The stochastic integral map is an isometry from $(\mc{B}, \norm{\cdot}_{2,r_1}')$ to  $L^2(\mc{A}, \phi)$.
Therefore the stochastic integral map can be continuously extended to the space  ${\mc{B}^{2, r_1}}'$.
\end{Prop}

\begin{proof}
Decompose $X(t) = r_1 t + X'(t)$, with $\phi[X'(t)] = 0$.
Let $N = \int V(s) \sharp dX(s)$, $M = \int U(s) \sharp dX(s)$.
Take $V = A \otimes B \chf_{[a, b)}$, $U = C \otimes D \chf_{[c, d)}$. Then 
\begin{align*}
\phi[N M^\ast] &= \phi[A (X'(b) - X'(a)) B D^\ast (X'(d) - X'(c)) C^\ast] \\
&\qquad + \phi[r_1 (\int V(s) ds) D^\ast (X'(d) - X'(c)) C^\ast] \\
&\qquad + \phi[A (X'(b) - X'(a)) B r_1 (\int U(s) ds)^\ast] + \phi[r_1 (\int V(s) ds) r_1 (\int U(s) ds)^\ast] \\
&= \abs{[a,b) \cap [c,d)} \phi[C^\ast A] \phi[B D^\ast] + 0 + 0 + r_1^2 \phi[(\int V(s) ds) (\int U(s) ds)^\ast] \\
&= \int \langle V(s), U(s) \rangle ds + r_1^2 \phi[(\int V(s) ds) (\int U(s) ds)^\ast].
\end{align*}
Both sides are linear in $U, V$, so the equality holds for arbitrary simple adapted biprocesses.
\end{proof}

\begin{Lemma}
\label{lem:tr}
For $U \in {\mc{B}^{2, r_1}}'$ such that $\sup_{0 \leq s \leq t} \norm{m(U)} < \infty$, 
\begin{equation}
\label{tr}
\phi[\int_0^t U(s) \sharp dX(s)] = \int_0^t \phi[m(U(s))] r_1 ds,
\end{equation}
where $m : \mc{A} \otimes \mc{A}^{op} \rightarrow \mc{A}$ is the usual multiplication.
\end{Lemma}
\begin{proof}
The formula holds for simple adapted biprocesses. Denote $M = \int_0^t U(s) \sharp dX(s)$, $N = r_1 \int_0^t m(U(s)) ds$. $\abs{\phi[M]} \leq \phi[M M^\ast] = \norm{M}_2^2$, so $\phi$ is continuous on $L^2(\mc{A}, \phi)$, and the left-hand-side of the equation \eqref{tr} is continuous. If $r_1 = 0$, both sides of the equation are $0$. If $r_1 \neq 0$, $r_1 \abs{\phi[N]} \leq \phi[N N^\ast] \leq \norm{U}_{2, r_1}^2$. Finally, since $\phi$ is normal, the hypotheses imply that $\phi[N] = r_1 \int_0^t \phi[m(U(s))] ds$.
\end{proof}

\section{It\^{o} formulas}
\label{sec:ito}
\begin{Prop}[It\^{o} product formula]
\label{prop:prIto}
Let $\set{U_j}_{j=1}^J$ and $\set{V_i}_{i=1}^I$ be two finite collections of biprocesses in $\mc{B}$.
Let $N = \sum_{i=1}^I N_i$, $M = \sum_{j=1}^J M_j$, where
\begin{equation*}
N_i = \int V_i(t) \sharp d \Delta_i(t), \qquad M_j = \int U_j(t) \sharp d \Delta_j(t).
\end{equation*}
Then 
\begin{align*}
N M &= \sum_{j=1}^J \int (N(t) \otimes 1) \sharp U_j(t) \sharp d \Delta_j(t) \\
&+ \sum_{i=1}^I \int (1 \otimes M(t)) \sharp V_i(t) \sharp d \Delta_i(t) \notag \\
&+ \sum_{m=2}^{I+J} \sum_{i, j = 1}^{i + j = m} \int (V_i(t) \otimes_2 U_j(t)) \sharp d \Delta_m(t). \notag
\end{align*}
\end{Prop}
\begin{proof}
It suffices to show that for  $N = \int V_i(t) \sharp d \Delta_i(t)$, $M = \int U_j(t) \sharp d \Delta_j(t)$,
\begin{multline}
\label{prIto2}
N M = \int (N(t) \otimes 1) \sharp U_j(t) \sharp d \Delta_j(t) \\
+ \int (1 \otimes M(t)) \sharp V_i(t) \sharp d \Delta_i(t) + \int (V_i(t) \otimes_2 U_j(t)) \sharp d \Delta_{i+j}(t).
\end{multline}
Moreover, we may assume that $V_j = A \otimes B \chf_{[0,a)}$ and either $U_i = C \otimes D \chf_{[a, b)}$ or $U = C \otimes D \chf_{[0,a)}$. In the first case it is easy to see that both sides of \eqref{prIto2} are equal to $A \Delta_i(a) B C (\Delta_j(b) - \Delta_j(a))D$.

In the second case, the statement to be proven is
\begin{multline*}
A \Delta_i(a) B C \Delta_j(a) D = \int (A \Delta_i(t) B C \otimes D \chf_{[0,a)}(t)) \sharp d \Delta_j(t) \\
+ \int ( A \otimes B C \Delta_j(t) D \chf_{[0, a)}(t)) \sharp d \Delta_i(t) + \int (A \otimes D \phi[B C] \chf_{[0, a)}(t)) \sharp d \Delta_{i+j}(t).
\end{multline*}
Now
\begin{equation*}
A \Delta_i(a) B C \Delta_j(a) D  = \lim_{N \rightarrow \infty} \sum_{k,l = 1}^N A X_k^i B C X_l^j D,
\end{equation*}
where the limit is taken in the operator norm, and $X_k$ ($= X_k^{(N)}$) is as in the definition of $\Delta$. The sum is equal to 
\begin{align*}
& \sum_{l=1}^N \sum_{k=1}^{l-1} A X_k^i B C X_l^j D + \sum_{k=1}^N \sum_{l=1}^{k-1} A X_k^i B C X_l^j D + \sum_{k=1}^N A X_k^i B C X_k^j D \\
&= \sum_{k=1}^N A \Delta_i \left( a\frac{k-1}{N} \right) B C X_k^j D + \sum_{k=1}^N A X_k^i B C \Delta_j \left( a\frac{k-1}{n} \right) + \sum_{k=1}^N A X_k^i B C X_k^j D.
\end{align*}
The operator norm limit of the third term is $A \phi[BC] \Delta_{i+j}(a) D$ by Corollary 14 of \cite{Ans00}. Now consider the first term. It can be written as $\int W_N(s) d\Delta_j(s)$, with
\begin{equation*}
W_N(s)  = \sum_{k=1}^N A \Delta_i \left( a\frac{k-1}{N} \right) B C \otimes D \chf_{[a\frac{k-1}{N}, a\frac{k}{N})}(s).
\end{equation*}
For $s \in [a\frac{k-1}{N}, a\frac{k}{N})$, by the proof of Lemma \ref{lem:norm} $\norm{\Delta_i(s) - \Delta_i(a\frac{k-1}{N})}_2 = O(\frac{1}{\sqrt{N}})$ uniformly in $k$. Therefore $W_N$ converges to $(A \Delta_i BC \otimes D \chf_{[0,a)})$ in the $2$-norm. This implies the convergence of integrals. The argument for the second term is similar.
\end{proof}

\begin{Thm}
\label{thm:Ito}
Let $M = \sum_{m=1}^K \int U_m \sharp d \Delta_m$ with $\set{U_i} \subset \mc{B}$, and $p$ be a polynomial. Then 
\begin{equation}
\label{ito1}
p \left( M \right) = \sum_{m=1}^\infty \sum_{k=1}^m \int \frac{1}{k!} \phi_{k+1} [\partial^k(p)(M) \sharp S(m,k)] \sharp d \Delta_m,
\end{equation}
where
\begin{equation*}
\sharp S(m,k) = \sum_{i(1), i(2), \ldots, i(k) = 1}^{i(1) + i(2) + \ldots + i(k) = m} \sharp \m_k(U_{i(1)}, U_{i(2)}, \ldots, U_{i(k)}).
\end{equation*}
\end{Thm}
\begin{proof}
It suffices to prove the formula for a monomial $p(x) = x^n$. Note that in this case the only non-zero terms on the right-hand-side of \eqref{ito1} are for $k \leq n$, $m \leq nK$, so the sum has finitely many terms. We proceed by induction on $n$. Denote by $U_{m, n}$ the coefficients in the expansion
\begin{equation*}
\left( \sum_{m=1}^K \int U_m \sharp d \Delta_m \right)^n = \sum_{m=1}^{nK} \int U_{m, n} \sharp d \Delta_m.
\end{equation*}
Using Proposition \ref{prop:prIto} with $N = M^n$, we know that 
\begin{equation}
\label{ind}
U_{m, n+1} = (M^n \otimes 1) \sharp U_m + (1 \otimes M) \sharp U_{m,n} + \sum_{i, j = 1}^{i+j=m} U_{i,n} \otimes_2 U_j,
\end{equation}
and we need to show that 
\begin{equation*}
U_{m, n} = \sum_{k=1}^m \frac{1}{k!} \phi_{k+1} [\partial^k M^n \sharp S(m,k)]
\end{equation*}
satisfy these equations. Indeed,

\begin{align*}
\sum_{k=1}^m \frac{1}{k!} \phi_{k+1} [\partial^k M^{n+1} \sharp S(m,k)] 
&= \sum_{k=1}^m \frac{1}{k!} \phi_{k+1} \left[\left( (1 \otimes \ldots \otimes M) \partial^k M^n + k \partial^{k-1} M^n \otimes 1 \right) \sharp S(m,k) \right] \\
&= (1 \otimes M) \sharp U_{m,n} + \sum_{k=1}^m \frac{1}{(k-1)!} \phi_{k+1} [(\partial^{k-1} M^n \otimes 1) \sharp S(m,k)].
\end{align*}
The second term of this sum is equal to
\begin{align*}
& (M^n \otimes 1) \sharp S(m,1)+ \sum_{k=1}^{m-1} \frac{1}{k!} \phi_{k+2} [(\partial^k M^n \otimes 1) \sharp S(m,k+1)] \\
&= (M^n \otimes 1) \sharp U_m + \sum_{k=1}^{m-1} \frac{1}{k!} \phi_{k+2} \left[(\partial^k M^n \otimes 1) \sharp \sum_{i(1), \ldots, i(k+1) =1}^{i(1) + \ldots + i(k+1) = m} \m_{k+1}(U_{i(1)}, \cdots, U_{i(k+1)}) \right] \\
\intertext{for $m > 1$. The second term of this sum, in turn, is equal to}
&\sum_{k=1}^{m-1} \frac{1}{k!} \phi_{k+2} \left[\sum_{i=k, j = 1}^{i+j=m} (\partial^k M^n \otimes 1) \sharp \sum_{i(1), i(2), \ldots, i(k) =1}^{i(1) + i(2) + \ldots + i(k) = i} \m_{k+1}(U_{i(1)}, U_{i(2)}, \cdots, U_{i(k)}, U_j)  \right] \\
&= \sum_{i, j = 1} ^{i+j=m} \phi_3 \left[ \m_2 \left( \sum_{k=1}^i \frac{1}{k!} \phi_{k+1} \left[\partial^k M^n \sharp \sum_{i(1), i(2), \ldots, i(k) =1}^{i(1) + i(2) + \ldots + i(k) = i} \m_k(U_{i(1)}, U_{i(2)}, \cdots, U_{i(k)}) \right],  U_j \right) \right] \\
&= \sum_{i, j = 1} ^{i+j=m} \phi_3 \left[ \m_2 \left( \sum_{k=1}^i \frac{1}{k!} \phi_{k+1} \left[  \partial^k M^n  \sharp S(i, k) \right],  U_j \right) \right] \\
&= \sum_{i, j = 1} ^{i+j=m} \left( \sum_{k=1}^i \frac{1}{k!} \phi_{k+1} \left[ \partial^k M^n \sharp S(i, k) \right] \right) \otimes_2 U_j \\
&= \sum_{i, j = 1} ^{i+j=m} U_{i,n} \otimes_2 U_j.
\end{align*}
Putting these three terms together we obtain the formula \eqref{ind}.
\end{proof}

\begin{Cor}[Functional It\^o formula]
\label{cor:Ito}
For $U \in \mc{B}$, $M = \int U \sharp dX$, the formula simplifies to
\begin{equation*}
p \left( \int U\sharp dX\right) = \sum_{k=1}^\infty \int \frac{1}{k!} \phi_{k+1} [ \partial^k(p)(M) \sharp \m_k(U, U, \ldots, U) ] \sharp d \Delta_k.
\end{equation*}
\end{Cor}

\begin{Remark}
There is a similar formula for the independence-based probability involving the usual derivatives of $p$. It does not seem to appear in the standard probability textbooks, but cf. \cite{Par92}.
\end{Remark}

\begin{Remark}
Theorem \ref{thm:Ito} has the following heuristic interpretation. Denote $\Delta_k = (dX)^k$. For $M$ as in the theorem, denote $dM = \sum_{m=1}^\infty U_m \sharp (dX)^m$. Then 
\begin{equation*}
d(p(M)) = p(M + dM) - p(M) = \sum_{k=1}^\infty \frac{1}{k!} \partial^k(p)(M) \sharp \m_k(dM, dM, \ldots, dM),
\end{equation*}
where in this context ``$\sharp \m_k(dM, dM, \ldots, dM)$'' means ``put $dM$ in place of each $\otimes$ and multiply through''. The theorem them follows from Corollary 14 of \cite{Ans00}, which says in this language that
\begin{align*}
Z_0 (dX)^{i(1)} Z_1 (dX)^{i(2)} \cdots (dX)^{i(n)} Z_n &= Z_0 \phi[Z_1] \cdots \phi[Z_{n-1}] (dX)^{\sum_{j=1}^n i(j)} Z_k \\
&= \phi_{k+1}[Z_0 \otimes Z_1 \otimes \ldots \otimes Z_k] \sharp (dX)^{\sum_{j=1}^n i(j)}
\end{align*}
\end{Remark}

\begin{Prop}
\label{prop:prod}
Let $V_i \in \mc{B}$, $M_i = \int V_i \sharp dX$, for $1 \leq i \leq n$. Then
\begin{multline*}
\prod_{i=1}^n M_i = \sum_{k=1}^n \int \phi_{k+1}[\sum_{i(1), \ldots i(k) = 1}^n (M_1 M_2 \cdots M_{i(1) - 1} \otimes M_{i(1) + 1} \cdots M_{i(2) - 1} \otimes \cdots M_{i(k) + 1} \cdots M_n) \\
\sharp \m_k(V_{i(1)}, V_{i(2)}, \ldots, V_{i(k)})] \sharp d \Delta_k.
\end{multline*}
\end{Prop}
\begin{proof}
This is a slight generalization of Corollary \ref{cor:Ito}; the proof proceeds as in Theorem \ref{thm:Ito}.
\end{proof}

\begin{Cor}
\label{cor:Ito2}
Even more particularly,
\begin{equation*}
p (X(t)) = \sum_{k=1}^\infty \int_0^t \frac{1}{k!} \phi_{k+1} [\partial^k(p)(X(s))] \sharp d \Delta_k(s).
\end{equation*}
\end{Cor}
\begin{proof}
$\chf_{[0,t)} 1 \otimes 1$ is a simple biprocess.
\end{proof}

\section{Extension to $L^\infty$}
\label{sec:linf}
Let $U$ be a scalar-valued biprocess, $U = \chf_{[0, \eps)}$. Then $M = \int U \sharp dX = X(\eps)$ and so $\norm{M} = \norm{\mu_\eps}$. This does not go to $0$ in the operator norm as $\eps \rightarrow 0$. So the stochastic integral map, considered as a map into $L^\infty(\mc{A}, \phi)$, is in general not continuous in the norm $\norm{\norm{U}_\infty}_2$, unlike in the free Brownian motion case. On the other hand, this suggests the use of the norm $\norm{\norm{U}_\infty}_\infty$.

Throughout this section, assume that all free cumulants $r_i$ of $\mu$ are non-negative. Such measures include the semicircular distribution, the free Poisson distribution, and more generally any compactly supported free compound Poisson distribution for which all the moments of its generator distribution are non-negative. Note that this condition does not imply that the corresponding operators are positive.

\begin{Defn}
Let $f \in L^{\infty-}(\mf{R}, dx)$, where $L^{\infty-}(\mf{R}, dx) = \bigcap_{p \geq 1} L^p(\mf{R}, dx)$. 
Define, for $n$ even,
\begin{equation*}
\norm{f}_{n, \mu} = \sqrt[n]{\sum_{\pi \in \NC(n)} \prod_{B_i \in \pi} r_{\abs{B_i}} \norm{f}_{\abs{B_i}}^{\abs{B_i}}}
\end{equation*}
and
\begin{equation*}
\norm{f}_{\infty, \mu} = \limsup_{n \rightarrow \infty} \norm{f}_{n, \mu}.
\end{equation*}
Let $\mc{B}^{n, \mu}$, for $n$ even or $\infty$, be the completion of $\mc{B}$ with respect to the norm $\norm{\norm{\cdot}_\infty}_{n, \mu}$. 
\end{Defn}

\begin{Prop}
\label{prop:norm}
$\norm{\norm{\cdot}_\infty}_{n, \mu}$ is indeed a norm.
\end{Prop}

\begin{Thm}
\label{thm:li}
Let $U \in \mc{B}$, $M = \int U(s) \sharp dX(s)$. Then for even $n$
\begin{equation*}
\norm{\int U(s) \sharp dX(s)}_n \leq \norm{\norm{U}_\infty}_{n, \mu}.
\end{equation*}
Therefore the stochastic integral map can be continuously extended to a contraction from $\mc{B}^{n, \mu}$ to $L^n(\mc{A}, \phi)$, and from $\mc{B}^{\infty, \mu}$ to $L^\infty(\mc{A}, \phi)$. 
\end{Thm}

\begin{Remark}
$\norm{U}_{2, r_1}' \leq \norm{U}_{2, \mu}$. Therefore $\mc{B}^{2, \mu} \subset {\mc{B}^{2, r_1}}'$.
\end{Remark}

We first prove the Theorem in the scalar case, when the inequality is in fact an equality. Note that in this case the non-negativity of the free cumulants is not necessary.

\begin{Prop}
\label{prop:scalar}
Let $f$ be a simple positive function, and $M = \int f(s) dX(s)$. Then the distribution of $M$ is the unique probability measure $\nu$ determined by $r_k(\nu) = r_k(\mu) \norm{f}_k^k$. In particular, 
\begin{equation*}
\phi[M^n] = m_n(\nu) = \norm{f}_{n, \mu}.
\end{equation*}
\end{Prop}
\begin{proof}
Denote $M(t) = \int_0^t f(s) dX(s)$. Since $f$ is simple, $M = M(t)$ for $t$ large enough, so it suffices to prove the Proposition for finite $t$. That is, we need to show that the distribution of $M(t)$ is the unique probability measure $\nu_t$ determined by $r_k(\nu_t) = r_k(\mu) \norm{f \chf_{[0,t)}}_k^k$. 
For any finite $t$, $M(t)$ is well-defined by the results in subsection \ref{sec:l2}.
Corollary \ref{cor:Ito} gives
\begin{equation}
\label{B2}
\phi[M(t)^n] = \int_0^t \sum_{k=1}^n \sum_{i(1), \ldots, i(k) = 0}^{i(1) + \ldots + i(k) = n-k} (i(1) + 1) \phi[M(s)^{i(1)}] \cdots \phi[M(s)^{i(k)}] f(s)^k r_k ds .
\end{equation}
In terms of generating functions equation \eqref{B2} corresponds to the differential equation
\begin{equation}
\label{B3}
\partial_t  G_{M(t)}(z) + \partial_z G_{M(t)}(z) f(t) \Rtr(f(t) G_{M(t)}(z)) = 0,
\end{equation}
where $\partial_z, \partial_t$ are the usual partial derivatives, $G_A(z) = \phi[(z - A)^{-1}]$ is the Cauchy transform of the distribution of $A$, and $\Rtr$ is the $R$-transform of $\mu$. See more on this in subsection \ref{sec:res}.

But $G_{\nu_t}$ is a solution to equation \eqref{B3}, with the same initial conditions. Indeed, that equation can be obtained by the same method as the original one in \cite{Voi86}. Namely, start with the equation
\begin{equation*}
G_{\nu_t} \left( \frac{1}{z} + \Rtr_{\nu_t}(z) \right) = z,
\end{equation*}
where $\Rtr_{\nu_t}(z) = \sum_{k=1}^\infty r_k (\norm{f \chf_{[0,t)}}_k^k) z^{k-1}$. Differentiating it with respect to $t$, we obtain
\begin{equation}
\label{B4}
\partial_t G_{\nu_t} \left( \frac{1}{z} + \Rtr_{\nu_t}(z) \right) + \partial_z G_{\nu_t} \left( \frac{1}{z} + \Rtr_{\nu_t}(z) \right) \partial_t \Rtr_{\nu_t}(z) = 0.
\end{equation}
But $\partial_t \Rtr_{\nu_t}(z) = \sum_{k=1}^\infty r_k f(t)^k z^{k-1} = f(t) \Rtr(f(t)z)$. Substituting $G_{\nu_t}(z)$ for $z$ in equation \eqref{B4}, we obtain equation \eqref{B3}. A posteriori, $\nu_t$ is the distribution of $M(t)$ and therefore a positive measure. Since $M(t)$ is bounded, its free cumulants determine a unique probability distribution.
\end{proof}

\begin{Remark}
The origin of the above proposition is in the combinatorial formulas of \cite{Ans00,Ans01b}. Also, the following heuristic argument should de-mystify the result. We note that since $f$ is a scalar-valued function, $\int_0^t f(s) dX(s)$ is a process with freely independent increments. Also, since the $R$-transform of a sum of freely independent variables is a sum of the $R$-transforms, heuristically, the same should be true for the integral. Denote by $D_c$ the dilation operator, $D_c \nu(\cdot) = \nu(c^{-1} \cdot)$. Then by stationarity
\begin{align*}
\sum_{k=1}^\infty r_k(\int f(s) dX(s)) z^{k-1} 
&= \Rtr_{\int f(s) dX(s)}(z) = \int \Rtr_{f(s) dX(s)}(z) \\
&= \int \Rtr_{D_{f(s)} \mu}(z) = \int f(s) \Rtr_\mu(f(s)z) \\
&= \int \sum_{k=1}^\infty r_k f(s)^k z^{k-1} = \sum_{k=1}^\infty (r_k \norm{f}_k^k) z^{k-1}.
\end{align*}
\end{Remark}

\begin{Cor}
Whenever $\norm{f}_{\infty, \mu} < \infty$, the $\limsup$ that defines it is in fact a limit, which equals to $\norm{\int f(s) dX(s)}$.
\end{Cor}

\begin{proof}[Proof of Proposition \ref{prop:norm}]
$\norm{f}_{n, \mu}$ is an $n$-th root of a homogeneous polynomial of degree $n$, with positive coefficients, in various $p$-norms of $f$, with no constant term. Therefore it is homogeneous and $0$ only at $0$. 

Now we prove the triangle inequality.
First take two functions $f \geq g \geq 0$ in $L^{\infty-}$. Since all the free cumulants $r_k$ are non-negative, $r_k \norm{f}_k^k \geq r_k \norm{g}_k^k$, and so $\norm{\int f(s) dX(s)}_p \geq \norm{\int g(s) dX(s)}_p$. 

Now let $U, V$ be simple biprocesses, $f(s) = \norm{U(s)} + \norm{V(s)}$, $g(s) = \norm{U(s) + V(s)}$. Then
\begin{align*}
\norm{\norm{U + V}_\infty}_{p, \mu} &= \norm{\int \norm{U(s) + V(s)}_\infty dX(s)}_p \\
&\leq \norm{\int (\norm{U(s)}_\infty + \norm{V(s)}_\infty) dX(s)}_p \\
&= \norm{\int \norm{U(s)}_\infty dX(s) + \int \norm{V(s)}_\infty dX(s)}_p \\
&\leq \norm{\int \norm{U(s)}_\infty dX(s)}_p + \norm{\int \norm{V(s)}_\infty dX(s)}_p \\
&= \norm{\norm{U}_\infty}_{p, \mu} + \norm{\norm{V}_\infty}_{p, \mu}.
\end{align*}
By approximation, $\norm{\norm{\cdot}_{\infty}}_{p, \mu}$ is a norm on $\mc{B}^{p, \mu}$.
For $p = \infty$, the triangle inequality and homogeneity follow by the limiting procedure, and for a non-negative function $f$, $\norm{f}_{\infty, \mu} = \norm{\int f(s) dX(s)} = 0$ iff $f = 0$. 
\end{proof}

\begin{Lemma}
\label{lem:moments}
For $1 \leq i \leq n$, let $V_i \in \mc{B}$, $M_i = \int V_i(s) \sharp dX(s)$, $N_i = \int \norm{V_i(s)} dX(s)$. Then
\begin{equation*}
\abs{\phi[M_1 M_2 \ldots M_n]} \leq \phi[N_1 N_2 \ldots N_n].
\end{equation*}
\end{Lemma}
\begin{proof}
We do this by induction. 
For $n=1$, by Lemma \ref{lem:tr}
\begin{align*}
\abs{\phi[\int V(s) \sharp dX(s)]} &= r_1 \abs{\int \phi[m(V(s))] ds} \\
&\leq r_1 \int \norm{m(V(s))} ds \leq r_1 \int \norm{V(s)} ds = \phi[\int \norm{V(s)} dX].
\end{align*}

By Proposition \ref{prop:prod},
\begin{multline*}
\phi[\prod_{i=1}^n M_i] \\
= \sum_{k=1}^n \int \phi \Bigl[ m \Bigl( \phi_{k+1} \Bigl[ \sum_{i(1), \ldots i(k) = 1}^n (M_1 M_2 \cdots M_{i(1) - 1} \otimes M_{i(1) + 1} \cdots M_{i(2) - 1} \otimes \cdots M_{i(k) + 1} \cdots M_n)  \\
\sharp \m_k(V_{i(1)}, V_{i(2)}, \ldots, V_{i(k)}) \Bigr] \Bigr) \Bigr] r_k ds,
\end{multline*}
and so
\begin{align*}
\abs{\phi[\prod_{i=1}^n M_i]} &\leq \sum_{k=1}^n \int \sum_{i(1), \ldots i(k) = 1}^n \abs{\phi[M_1 M_2 \cdots M_{i(1) - 1} M_{i(k) + 1} \cdots M_n]} \\
& \qquad \times \abs{\phi[M_{i(1) + 1} \cdots M_{i(2) - 1}]} \ldots \abs{\phi[M_{i(k-1) + 1} \cdots M_{i(k) - 1}]} \prod_{j=1}^k \norm{V_{i(j)}} r_k ds \\
&\leq \sum_{k=1}^n \int \sum_{i(1), \ldots i(k) = 1}^n \phi[N_1 N_2 \cdots N_{i(1) - 1} N_{i(k) + 1} \cdots N_n] \\
& \qquad \times \phi[N_{i(1) + 1} \cdots N_{i(2) - 1}] \ldots \phi[N_{i(k-1) + 1} \cdots N_{i(k) - 1}] 
\prod_{j=1}^k \norm{V_{i(j)}} r_k ds \\
&= \phi[\prod_{i=1}^n N_i],
\end{align*}
where we have used the induction hypothesis as well as the equality for positive scalar functions.
\end{proof}

\begin{proof}[Proof of Theorem \ref{thm:li}]
Apply Lemma \ref{lem:moments} to $M_i = M$ for $i$ odd, $M_i = M^\ast$ for $i$ even. We get
\begin{equation*}
\abs{\phi[\abs{M}^{2n}]} \leq \phi[N^{2n}] = \norm{\norm{U}_\infty}_{2n, \mu},
\end{equation*}
where in the second equality we have used Proposition \ref{prop:scalar}.
\end{proof}

\begin{Ex}
If $\mu$ is the semicircular distribution, then $\norm{f}_{2n, \mu} = \sqrt[2n]{c_n \norm{f}_2^{2n}} = \sqrt[2n]{c_n} \norm{f}_2$ (where $c_n$'s are the Catalan numbers) and so $\norm{f}_{\infty, \mu} = 2 \norm{f}_2$. In this case $r_k(\nu) = \delta_{k2} \norm{f}_2^2$, and so $\nu = \mu_{\norm{f}_2^2}$. So $\norm{U dX} \leq 2 \norm{\norm{U}_\infty}_2$. So in this case we recover the result of \cite{BS98}, in fact with a slightly better constant, which by our results is optimal.
\end{Ex}

\begin{Ex}
If $f(t) = \chf_{[a,a + \eps)}$, then $r_k(\nu) = r_k \eps$, and so $\nu = \mu_\eps$, $\norm{f}_{\infty, \mu} = \norm{\mu_\eps}$. In particular, $\norm{\int_0^t U dX} \leq \norm{\norm{U}_\infty}_\infty \norm{\mu_t}$. 
\end{Ex}

\begin{Cor}
If $\norm{\norm{U}_\infty}_k < C$ for all $k \leq n$, then 
\begin{equation*}
\norm{\int U dX}_n \leq \norm{\norm{U}_\infty}_{n, \mu} \leq C m_n(\mu).
\end{equation*}
\end{Cor}

\subsection{More on It\^{o} formulas}
Next we extend the product and the functional It\^{o} formulas to biprocesses in $\mc{B}^{\infty, \mu}$.

The following is a very preliminary form of the lower bound in the Burkholder-Davis-Gundy inequality.
\begin{Lemma}
Denote by $\mu^k$ the distribution of $\Delta_k$. Then for simple $f$, 
\begin{equation*}
\norm{f}_{n, \mu^k} \leq \norm{\abs{f}^{1/k}}_{nk, \mu}^k.
\end{equation*}
\end{Lemma}
\begin{proof}
For $\nu^k$ the distribution of $\int f(s) d \Delta_k(s)$, $r_i(\nu^k) = r_{ik} \norm{f}_i^i$. Since all the free cumulants are positive, 
\begin{equation*}
\sum_{\pi \in \NC(n)} \prod_{B \in \pi} r_{\abs{B}k} \norm{f}_{\abs{B}}^{\abs{B}} \leq \sum_{\pi \in \NC(nk)} \prod_{B \in \pi} r_{\abs{B}} \norm{f^{1/k}}_{\abs{B}k}^{\abs{B}k}.
\end{equation*}
Therefore $\norm{f}_{n, \mu^k} \leq \norm{\abs{f}^{1/k}}_{nk, \mu}^k$. In particular $\norm{f}_{\infty, \mu^k} \leq \norm{\abs{f}^{1/k}}_{\infty, \mu}^k$.
\end{proof}

\begin{Prop}
Let $\set{U_j}_{j=1}^\infty$ and $\set{V_i}_{i=1}^\infty$ be two collections of biprocesses in $\mc{B}^{\infty, \mu}$, with $ \sum_{j=1}^\infty \norm{\norm{U_j}_\infty^{1/j}}_{\infty, \mu}^j < \infty$, $\sum_{i=1}^\infty \norm{\norm{V_i}_\infty^{1/i}}_{\infty, \mu}^i < \infty$.
Let $N = \sum_{i=1}^\infty N_i$, $M = \sum_{j=1}^\infty M_j$, where
\begin{equation*}
N_i = \int V_i(t) \sharp d \Delta_i(t), \qquad M_j = \int U_j(t) \sharp d \Delta_j(t).
\end{equation*}
Then 
\begin{align}
N M &= \sum_{j=1}^\infty \int (N(t) \otimes 1) \sharp U_j(t) \sharp d \Delta_j(t) \label{prIto} \\
&+ \sum_{i=1}^\infty \int (1 \otimes M(t)) \sharp V_i(t) \sharp d \Delta_i(t) \notag \\
&+ \sum_{m=2}^\infty \sum_{i, j = 1}^{i + j = m} \int (V_i(t) \otimes_2 U_j(t)) \sharp d \Delta_m(t). \notag
\end{align}
\end{Prop}
\begin{proof}
First assume $U_i, V_j \in \mc{B}$.
\begin{equation*}
\norm{\int U_j(t) \sharp d \Delta_j(t)} \leq \norm{\norm{U_j}_\infty^{1/j}}_{\infty, \mu}^j,
\end{equation*}
so by hypothesis the series defining $M$ and $N$ converge absolutely. Also,
\begin{align*}
\norm{\int (N_i(s) \otimes 1) \sharp U_j(t) \sharp d \Delta_j(t)} 
&\leq \norm{N_i} \norm{\int U_j(t) \sharp d \Delta_j(t)} \\
&\leq \norm{\norm{U_j}_\infty^{1/j}}_{\infty, \mu}^j \norm{\norm{V_i}_\infty^{1/i}}_{\infty, \mu}^i.
\end{align*}
Finally,
\begin{multline*}
N_i M_j = \int (N_i(t) \otimes 1) \sharp U_j(t) \sharp d \Delta_j(t) \\
+ \int (1 \otimes M_j(t)) \sharp V_i(t) \sharp d \Delta_i(t) + \int (V_i(t) \otimes_2 U_j(t)) \sharp d \Delta_{i+j}(t).
\end{multline*}
implies that 
\begin{multline*}
\norm{\int (V_i(t) \otimes_2 U_j(t)) \sharp d \Delta_{i+j}(t)} \\
\leq \norm{N_i M_j} + \norm{\int (N_i(t) \otimes 1) \sharp U_j(t) \sharp d \Delta_j(t)} + \norm{\int (1 \otimes M_j(t)) \sharp V_i(t) \sharp d \Delta_i(t)} \\
\leq 3 \norm{\norm{U_j}_\infty^{1/j}}_{\infty, \mu}^j \norm{\norm{V_i}_\infty^{1/i}}_{\infty, \mu}^i.
\end{multline*}
So the sums on the right-hand-side of the equation \eqref{prIto} converge absolutely as well. The formula holds for finitely many $U_j, V_i \in \mc{B}$, and since both sides of that equation are continuous, the formula can be extended to $\mc{B}^{\infty, \mu}$.
\end{proof}

\begin{Prop}
\label{prop:anIto}
Let $U \in \mc{B}^{\infty, \mu}$, $M = \int U(s) \sharp dX(s)$. Let $p(x) = \sum_{n=1}^\infty a_n x^n$, with the series absolutely convergent for $\abs{x} < R$, with $R > 2\norm{\norm{U}_\infty}_{\infty, \mu} $. Then
\begin{equation}
\label{fIto2}
p \left( \int U\sharp dX\right) = \sum_{k=1}^\infty \int \frac{1}{k!} \phi_{k+1} [ \partial^k(p)(M) \sharp \m_k(U, U, \ldots, U) ] \sharp d \Delta_k.
\end{equation}
\end{Prop}
\begin{proof}
Let $\partial^k(p) = \sum_{n=1}^\infty a_n \partial^k(x^n)$, taken a priory as a formal power series. Assume that $U \in \mc{B}$. 
\begin{align*}
\norm{\int \frac{1}{k!} \phi_{k+1} [\partial^k(x^n)(M(s)) \sharp \m(U(s))] \sharp d \Delta_k(s)} 
&\leq \binom{n}{k} \norm{\int \norm{M(s)}^{n-k} \norm{U(s)}^k d \Delta_k(s)} \\
&\leq \binom{n}{k} \norm{\norm{U}_\infty}_{\infty, \mu}^{n-k} \norm{\norm{U}_\infty}_{\infty, \mu}^k \\
&= \binom{n}{k} \norm{\norm{U}_\infty}_{\infty, \mu}^n.
\end{align*}
For $\norm{\norm{U}_\infty}_{\infty, \mu} < R/2$, 
\begin{equation*}
\sum_{n=0}^\infty \sum_{k=1}^\infty \abs{a_n} \binom{n}{k} \norm{\norm{U}_\infty}_{\infty, \mu}^n 
\leq \sum_{n=0}^\infty \abs{a_n} (2 \norm{\norm{U}_\infty}_{\infty, \mu})^n < \infty.
\end{equation*}
Therefore the sum on the right-hand-side of equation \eqref{fIto2} converges absolutely. By continuity, the formula holds also for $U \in \mc{B}^{\infty, \mu}$.
\end{proof}

\subsection{Resolvent}
\label{sec:res}

In particular, the function $f(x) = (z - x)^{-1} = \sum_{k=0}^\infty z^{-(k+1)} x^k$ is analytic in $x$ for $\abs{x} < \abs{z}$. 
Denote by $\Res$ the resolvent function,
\begin{equation*}
\Res_{A}(z) = (z - A)^{-1}.
\end{equation*}
It has the following nice behavior with respect to the derivation $\partial$:
\begin{equation*}
\partial \Res_{A}(z) = \Res_{A}(z) \otimes \Res_{A}(z)
\end{equation*}
and more generally
\begin{equation*}
\partial^k \Res_{A}(z) = k! \Res_{A}(z)^{\otimes(k+1)}.
\end{equation*}
See \cite{Voi00a} for much deeper analysis of this property. For the resolvent, Proposition \ref{prop:anIto} and the appropriate modification of Corollary \ref{cor:Ito2} read, respectively,

\begin{equation*}
\Res_{M(t)}(z) = \sum_{k=1}^m \int_0^t \phi_{k+1} [\Res_{M(s)}(z)^{\otimes(k+1)} \sharp \m(U, U, \ldots, U)] \sharp d \Delta_k(s)
\end{equation*}
and

\begin{align}
\label{Burger}
\Res_{X(t)}(z) &= \sum_{k=1}^\infty \int_0^t \phi_{k+1} [\Res_{X(s)}(z)^{\otimes(k+1)}] \sharp d \Delta_k(s) \notag \\
&= \sum_{k=1}^\infty \int_0^t G_{X(s)}(z)^{k-1} \Res_{X(s)}(z) d \Delta_k(s) \Res_{X(s)}(z).
\end{align}
Here $G_{X(t)}(z) = \phi[\Res_{X(t)}(z)]$ is the Cauchy transform of the distribution $\mu_t$ of $X(t)$.

Evaluating $\phi$ on both sides of \eqref{Burger}, we obtain by Lemma \ref{lem:r_n}
\begin{align}
\label{ql}
0 &= G_{X(t)}(z) + \sum_{k=1}^\infty r_k \int_0^t G_{X(s)}(z)^{k-1} \partial_z G_{X(s)}(z) ds \notag\\
&= G_{X(t)}(z) + \int_0^t \Rtr(G_{X(s)}(z)) \partial_z G_{X(s)}(z) ds.
\end{align}
This is the well-known quasi-linear equation for $G_{X(t)}$ in an integrated form.

\begin{Remark}
This is a manifestation of a more general fact. It seems that a number of the combinatorial properties of free independence can be expressed on the level of operators, with partition-dependent stochastic measures of \cite{Ans00} replacing the scalar-valued $R$-transforms. The scalar identities are then obtained by taking the expectations with respect to the trace $\phi$. We will say more about this topic elsewhere.
\end{Remark}

\providecommand{\bysame}{\leavevmode\hbox to3em{\hrulefill}\thinspace}


\begin{thebibliography}{BNTr00}

\bibitem[Ans00]{Ans00}
Michael Anshelevich, \emph{Free stochastic measures via noncrossing
  partitions}, Adv. Math. \textbf{155} (2000), no.~1, 154--179,
  arXiv:math.OA/9903084.

\bibitem[Ans01]{Ans01b}
Michael Anshelevich, \emph{Free stochastic measures via noncrossing partitions
  {II}}, preprint, 2001, arXiv:math.OA/0102062.

\bibitem[BNTr00]{BNT00}
O.E. Barndorff-Nielsen and S.~Thorbj{\o}rnsen, \emph{Selfdecomposability and
  {L}\'{e}vy processes in free probability}, Odense preprint no. 24, 2000.

\bibitem[BS98]{BS98}
Philippe Biane and Roland Speicher, \emph{Stochastic calculus with respect to
  free {B}rownian motion and analysis on {W}igner space}, Probab. Theory
  Related Fields \textbf{112} (1998), no.~3, 373--409.

\bibitem[CD99]{CD99}
Thierry Cabanal-Duvillard, \emph{Fluctuations de la loi empirique de grandes
  matrices al\'{e}atoires}, preprint, 1999.

\bibitem[CDG00]{CDG00}
T.~Cabanal-Duvillard and A.~Guionnet, \emph{Large deviations upper bounds and
  non commutative entropies for some matrices ensembles}, ENS DMA preprint
  00-03, 2000.

\bibitem[Par92]{Par92}
Kalyanapuram~Rangachari Parthasarathy, \emph{A quantum stochastic approach to
  {I}t\^o's formula for {L}\'evy processes}, C. R. Acad. Sci. Paris S\'er. I
  Math. \textbf{315} (1992), no.~13, 1417--1420.

\bibitem[RW97]{RW97}
Gian-Carlo Rota and Timothy~C. Wallstrom, \emph{Stochastic integrals: a
  combinatorial approach}, Ann. Probab. \textbf{25} (1997), no.~3, 1257--1283.

\bibitem[Voi86]{Voi86}
Dan Voiculescu, \emph{Addition of certain noncommuting random variables}, J.
  Funct. Anal. \textbf{66} (1986), no.~3, 323--346.

\bibitem[VDN92]{VDN92}
D.~V. Voiculescu, K.~J. Dykema, and A.~Nica, \emph{Free random variables}, CRM
  Monograph Series, vol.~1, American Mathematical Society, Providence, RI,
  1992, A noncommutative probability approach to free products with
  applications to random matrices, operator algebras and harmonic analysis on
  free groups.

\bibitem[Voi00a]{Voi00a}
Dan Voiculescu, \emph{The coalgebra of the free difference quotient and free
  probability}, Internat. Math. Res. Notices \textbf{2000}, no.~2, 79--106.

\bibitem[Voi00b]{Voi00b}
Dan Voiculescu, \emph{Lectures on free probability theory}, Lectures on
  probability theory and statistics (Saint-Flour, 1998), Springer, Berlin,
  2000, pp.~279--349.

\end{thebibliography}

\end{document}